\documentstyle{amsppt}
\input amsppt1
\nologo
%\NoPageNumbers
%
% ------ Macros ------
%
\font\b=cmr10 scaled \magstep4
\def\bigzerou{\smash{\kern-25pt\lower1.7ex\hbox{\b 0}}}
\hsize=350pt
\vsize=544pt
\hbadness=5000
\tolerance=1000
\NoRunningHeads
\magnification=\magstep1
%\pageno=27
\def\no{\noindent}
\def\f{\frak}

\def\Q{\Bbb Q}

\def\al{\alpha}

\def\F5{\Bbb F_5}

\def\Z{\Bbb Z}

\def\m{\frak m}

\def\fp{\frak p}
\def\ff{\frak f}
\def\l{\frak l}

\def\ml{\frak m\frak l}

\def\i{\infty}

\def\ga{\gamma}

\def\ot{\otimes}
\def\ov{\overline}

\topmatter

%\pageno=27
\topmatter
\title Euler systems of $K_2$ of CM elliptic curves \endtitle
\author
Kenichiro Kimura
\endauthor
 
\abstract
{We construct certain systems of elements in $K_2$ 
of some CM elliptic curves. When the classnumber of the 
field of CM is 1, The image of this system under the regulator
map forms an euler system in the sense of Rubin \cite{Ru}.}

\endtopmatter

\vskip -3ex
\noindent\hskip 2em
$\ssize{\text{
}}$
\vskip 3ex
\CenteredTagsOnSplits

% ------  Document  ------
%
\document
\baselineskip=13pt 
\subhead \S 1. Systems of $K_2$ of CM elliptic curves\endsubhead

Let $K$ be an imaginary quadratic field.
 Take a Hecke character $\phi$ of $K$ of type 
(1,0) and let $\ff_\phi$ be its conductor.

\no Let $H=K(\ff_\phi)$ be the ray class field of $K$ 
modulo $\ff_\phi$. Then by (\cite{dS}, lemma 1.4, Ch. 2)
there is an elliptic curve $E$ defined over $H$
with complex multiplication by the ring of integers
$\Cal O_K$ of $K$ such that
the Hecke character $\psi$ associated to $E$ is of the form

\hfil$\psi=\phi\circ N_{H/K}.$\hfil

\no We assume that $E(\Bbb C)$ is isomorphic to $\Bbb C/\Cal O_K$.

\no Let $\frak F$ be the conductor of $E$ and write $\f f=N_{H/K}(\f F)$.

\no For an integral ideal $\f m$ of 
$\Cal O_K$ $E[\f m]$ denotes the group of $\f m$-torsion points of 
$E(\ov{H})$. 

\no For any ideal $\f n$ of $\Cal O_K$ we write $H(\f n)$ for 
$H(E[\f n])$.

\no Let $p>3$ be a rational prime which splits in $K/\Q$, say
$p=\fp\bar{\fp}$. We assume that $p$ 
is relatively prime to $\frak F$.

\no Let $a$ 
be an integer relatively prime to $p$. Take a function 
$g_a\in H(E)$ such that div$g_a=a^2(0)-E[a]$. 

\no Let $\Cal L$ be the set of principal prime ideals of 
$\Cal O_K$ relatively prime to $\f f \bar{\fp} a$ and which
split completely in $H$.

\no Let $R$ be the set of ideals of $\Cal O_K$ which are divisible
only by primes in $\Cal L$.

%\no Let $\psi: (\text{ideal group of}\,\,H) \to K^*$ 
%be the Hecke character associated to 
%$E$. Then there exists a Hecke
%character $\phi$ of $K$ of type (1,0) such that 
%$\psi=\phi\circ N_{H/K}$er $\phi$. 

\no For each $\f m\in R$ we fix a generator 
$x_\f m$ of $E[\f m]$ as an $\Cal O_K$ module so that they satisfy the relation 
$[\phi(\f n)]x_{\f m\f n}=x_{\f m}$ for any $\f n\in R$.

\no Here $[\phi(\f n)]$ is $\phi(\f n)$ regarded as an element 
of $End(E)$. 

\no  Fix a generator $x_{\f f}$ of $E[\f f]$
as an $\Cal O_K$ module.
Since any ideal $\m\in R$ is prime to $\ff$, 
$[\phi(\m)]$ is an automorphism of $E[\f f]$.

\no For each $\f m\in R$ let $y_{\m}:=x_{\m}+[\phi(\m)]^{-1}x_{\f f}\in 
E(\ov{H})$  and take a function

\no $s_{\f m}\in H(\f m\ff)(E)$(the function field 
of $E\ot_HH(\f m\ff)$) such that

\no div$\,s_{\f m}=N(\m\f f)(y_{\f m})-N(\m\f f)(0).$ 

\no  For each $\ga \in E[a]-\{0\}$ take a function $t_\ga \in H(a)(E)$
 such that div $t_{\ga}=a(\ga)-a(0)$.

\no For any ideal $\m\in R$ consider the element 
$$\align&\al'_{\f m}:=N(\m\f f)^{-1}
\biggr(a\{g_a(y_{\f m})^{-1}g_a,\,\,s_{\f m}\}
-\sum_{\ga \in E[a]-\{0\}}\{s_{\f m}(\ga),\,t_{\ga}\}\biggl)\\
&\in \Gamma({E_{H(\f ma\ff)}}_{zar}
,\, \Cal K_2)\ot \Q\endalign$$ where $\Cal K_2$ is the Zariski sheaf
associated to presheaf
$$E\supset_{\text{open}} U\mapsto K_2(U)$$
and let 
$$ \al_{\f m}:=
\cases &[\phi(\f m)]_*N_{H(\m a\ff)/H(\m)}\al'_{\f m}\quad \fp |\m\\
       &[\phi(\f m\fp)]_*N_{H(\m\fp a\ff)/H(\m)}
         \al'_{\f m\fp}\quad \fp \nmid\m.
\endcases$$ 
Noting that 
$K_2$ of number fields is torsion, it can be checked
that the definition of $\al_{\m}$ is independent of choice of $g_a$, $s_{\m}$
and $t_{\ga}$. This definition
is similar to the one in Chap. 7 of \cite{Bl-Ka}.

We will state the main result.

\proclaim{Theorem} Let  $\m\in R$ and let $\l\in \Cal L$.
Then the element $\al_{\m\l}$ satisfies
the following equality:

$$\align    &(E1)\quad N_{H(\f m\f l)/H(\f m)}\al_{\f m\f l}=\al_{\f m}
\quad\text{if}\quad \f l|\f m \,\,\,\text{or}\,\,\,\l=\fp\\  
 &(E2)\quad N_{H(\f m\f l)/H(\f m)}
\alpha_{\f m\f l}=(1-[\phi(\f l)]_*Fr_{\f l}^{-1})\alpha_{\f m}
\quad\text{if}\quad\f l\nmid \f m\fp   .\endalign$$
Here $Fr_\l$ is the Frobenius element of any prime of $H$ 
over $\f l$.
\endproclaim

\demo{proof}

% Since the restriction
%$\Gamma({E_{H(\f m)}}_{zar}
%,\, \Cal K_2)\ot \Q \to \Gamma(U_{zar}
%,\, \Cal K_2)\ot \Q$ is injective, it suffices to show the relations (E1) and 
%(E2) in the latter space.

First we prove (E1). When $\l =\fp$ and 
$\fp\nmid \m$ this holds by definition.
So we assume that $\l | \m$. We will show the equality
$$[\phi(\l)]_*(N_{H(\ml a\ff)/H(\m a\ff)}\al_{\ml}')=\al_{\m}'.$$

\no Take functions $s_{\m\l}\in H(\ml \ff)(E)^*\ot \Q$ and 
$s_{\m}\in H(\m \ff)(E)^*\ot \Q$ such that 

\no div $s_{\m\l}=(y_{\m\l})-(0)$, div $s_{\m}=(y_{\m})-(0)$
respectively and such that $[\phi(\l)]_*s_{\m\l}=s_{\m}$. 

\no Since $Gal(H(\m\l a\ff)/H(\m a\ff))\simeq \Cal O_K/\l$ we see that
$$[\phi(\l)]^{-1}(y_{\m})=\underset{\tau\in Gal(H(\m\l a\ff)/H(\m a \ff))}
\to {\bigcup}y_{\m\l}^\tau.$$
Take a function $g_{\l}\in H(\m \ff)(E)^*\ot\Q$ such that
$ N_{H(\ml\ff)/H(\m\ff)}s_{\ml}=[\phi(\l)]^*(s_{\m})g_{\l}.$ Note that 
div $g_{\l}=\underset{c\in E[\l]}\to{\sum}(c)-N(\l)(0).$
 Then we have the equality

$$\align &N_{H(\ml a\ff)/H(\m a\ff)}\al'_{\ml}\\
&=a\{g_a,\,N_{H(\ml a\ff)/H(\m a\ff)}s_{\ml}\}
-\sum_{\ga\in E[a]-\{0\}}\{N_{H(\ml a\ff)/H(\m a\ff)}s_{\ml}(\ga),
\,t_{\ga}\}\\
&+\sum_{\tau\in Gal(H(\ml)/H(\m))}a\{g^{-1}_a(y_{\ml}^\tau),\,s_{\ml}^\tau\}\\
&= a\{g_a,\,[\phi(\l)]^*(s_{\m})g_{\l}\}
-\sum_{\ga\in E[a]-\{0\}}\{\bigr([\phi(\l)]^*(s_{\m})
g_{\l}\bigl)(\ga),\,t_{\ga}\}\\
&+\sum_{\tau\in Gal(H(\ml a\ff)/H(\m a\ff))}
a\{g^{-1}_a(y_{\ml}^\tau),\,s_{\ml}^\tau\}.
\endalign$$

\no Since $[\phi(\l)]_*s_{\ml}^\tau =s_{\m}^\tau =s_{\m}$ for 
$\tau\in  Gal(H(\ml a\ff)/H(\m a\ff))$, the equality

$$\align&[\phi(\l)]_*\sum_{\tau\in Gal(H(\m\l a\ff)/H(\m a\ff))}
\{g_a(y_{\m\l}^\tau)^{-1},\,s_{\m\l}^\tau\}\\
&=\{\prod_{\tau\in Gal(H(\m\l a\ff)/H(\m a\ff
)}g_a(y_{\m\f l}^\tau)^{-1},\,s_{\m}\}\\
&=\{[\phi(\f l)]_*g_a(y_{\m})^{-1},\,s_{\m}\}\endalign$$ holds.
Here the first equality holds by the projection formula.

$\quad$

\no Let $c\in E[\l]$ be a nonzero point 
and take a function $u\in H(\l)(E)^*\ot\Q$ which has 
the divisor $(c)-(0)$. Let 
$$A:=a\{g_a,\,g_{\l}\}
-\sum_{\ga\in E[a]-\{0\}}\{g_{\l}(\ga),
\,t_{\ga}\}\,\,\,\text{ and }\,\,\,
B:=\sum_{\xi\in Gal(H(\l)/H)}a\{g_a(c^\xi),\,u^\xi\}$$ be the elements of 
$K_2(H(\m\l a\ff)(E))\ot \Q.$ It can be seen that 
$A-B\in \Gamma (E_{H(\ml a\ff)},\,\Cal K_2)\ot \Q.$ We know that $[-1]$ 
acts on this group by $-1$. 
However, since $[-1]^*t_{\ga}/t_{-\ga}\in H(a)$, 
$[-1]^*g_a=\pm g_a$ and $[-1]^*g_{\f l}=
\pm g_{\f l}$, $A-B$ is invariant under $[-1]$.
From this we see that $A=B.$
Since $B\in Ker([\phi(\l)]_*)$, it follows that 
$A\in Ker([\phi(\l)]_*).$ 

\no Using this fact and the projection formula we get the equality $(E1)$.

Next we prove (E2). We will show that
$$[\phi(\l)]_*(N_{H(\ml a\ff)/H(\m a\ff)}\al_{\ml}'+Fr_{\l}^{-1}(\al_{\m}'))
=\al_{\m}'.$$
Since $(\m,\,\,\l)=1, \,\,\,[\phi(\l)]$ is an 
automorphism of $E[\m\ff]$.  Let $[\phi(\l)]^{-1}y_{\m}=:n\in E[\m]$ 
and take $s_{n}\in H(\m\ff)(E)^*\ot\Q$ such that
div$s_{n}=(n)-(0)$ and $[\phi(\l)]_*s_{n}=s_{\m}.$ Then the equality

$$Fr_{\l}^{-1}\al'_{\m}=a\{g_a(n)^{-1}g_a,\,\,s_{n}\}
-\sum_{\ga\in E[a]-\{0\}}\{s_{n}(\ga),\,t_{\ga}\}$$
holds.

\no Since $Gal(H(\m\l a\ff)/H(\m a\ff))\simeq (\Cal O_K/\l)^*$,
we see that
$$[\phi(\l)]^{-1}(y_\m)=\underset{\tau\in Gal(H(\m\l a\ff)/H(\m a\ff))}
\to{\bigcup}y_{\ml}^\tau \cup\{n\}.$$
Similarly as in the case
(E1) we take $s_{\ml}$ such that $[\phi(\l)]_*s_{\ml}=s_{\m}$ and 
let $g_{\l}\in H(\m\ff)(E)^*\ot\Q$ be the function satisfying the relation
$$N_{H(\ml \ff)/H(\m \ff)}(s_{\ml})s_{n}=[\phi(\l)]^*(s_{\m})g_{\l}.$$
The rest of the proof is similar to that of (E1).\qed\enddemo

\subhead{\S 2. Images under the regulator map}\endsubhead

Let $T_p(E)$ be the Tate module of $E$. There is a decomposition

\hfil $T_p(E)=T_{\fp}(E)\oplus T_{\bar{\fp}}(E)$ \hfil

\no where 

$$T_{\fp}(E)=\varprojlim_{n}E[\fp^n]$$

\no and

$$T_{\bar{\fp}}(E)=\varprojlim_{n}E[\bar{\fp}^n].$$

\no We will define a map 
$$r_G:\quad \Gamma(E_{H(\m)},\,\Cal K_2)\to H^1(H(\m)
,\,T_{\fp}(E)(1))$$ for each $\m\in R.$

\no Soul\'e defined the Chern class map
$$K_2(E_{H(\m)})\to H^2(E_{H(\m)},\,\Z/p^n(2))$$ in \cite{So}. 

\no Since $p$ is relatively prime to $\ff$, the prime  $\bar{\fp}$ of $K$
 is 
unramified in $H(\m)$ so that the field $H(\m)$ and the cyclotomic 
field $K(\mu_{p^n})$ is linearly disjoint over $K$.

\no Hence the group $H^0(H(\m),\,H^2(E_{\ov{H}},\,\Z_p/p^n(2)))=0.$

\no By Hochschild-Serre spectral sequence there is a map
$$Ch_n:\,\,K_2(E_{H(\m)})\to H^1(H(\m),\,H^1(E_{\ov{H}},\,\Z_p/p^n(2)))$$
and taking projective limit we get the map
$$\align Ch:\,\,K_2(E_{H(\m)})&\to \underset{n}\to{\varprojlim}
H^1(H(\m),\,H^1(E_{\ov{H}},\,\Z_p/p^n(2)))\\
&=H^1(H(\m),\,H^1(E_{\ov{H}},\,\Z_p(2))).
\endalign$$
Here the last equality follows from (\cite{Ta}, Proposition 2.2).

\no Take an element $\al\in \Gamma(E_{H(\m)} ,\,\Cal K_2)$. Since

\hfil$K_2(H(\m)(E))=\underset{E\supset_{\text{open}}U
}\to{\varinjlim}K_2(U)$\hfil

\no there is an open
set $U$ of $E$ and an element $\al_U\in K_2(U)$ such that
the image of $\al_U$ in $ \Gamma(E_{H(\m)} ,\,\Cal K_2)$ is $\al$.

\no Since $K_2$ of number fields is torsion, 
$\al_U$ is well defined modulo torsion elements. Noting that
there is an exact sequence
$$\underset{t\in E-U}\to{\oplus}K_2(\kappa(t))\to K_2(E)\to K_2(U)
\to \underset{t\in E-U}\to{\oplus}
K_1(\kappa(t))$$ and that $\al$ is in the 
kernel of Tame symbol, we see that there is an element
$\tilde{\al}\in K_2(E_{H(\m)} )$ such that $\tilde{\al}|_U=\al_U$.

\no Since the group
$H^1(H(\m),\,T_{\fp}(E)(1))$ has no torsion we can define
$r_G(\al)$ to be the projection of $Ch(\tilde{\al})$ to 
$H^1(H(\m),\,T_{\fp}(E)(1))$.

\no Since $\phi(\m)/N(\m)\in (\Cal O_K)_{\fp}^*$ for each $\m\in R$,
 $\al_{\m}$ defines a class $r_G(\al_{\m})
\in H^1(H(\m),\,T_{\fp}(E)(1)).$

\no When the class number of $K=1$, 
it can be checked that this family of cohomology classes forms 
an ``euler system'' in the sense of (\cite{Ru}, definition 2.1.1).

\no We let the field $\Cal K$ resp. $K_\i$ in loc.cit be the union of fields
$K(E[\m])$ for $\m\in R$ resp. the maximal $\Z_p$ extension
of $K$ in $K(E[\fp^\i])$.

\no The condition (i) in loc.cit can be checked using Proposition 
1.6 in Chap. II of \cite{dS} and the fact that the action of inertia group
of a prime of $K$ on $E[\f q]$ for any prime $\f q \nmid \ff$ factors through
the roots of unity in $K$.

\no The second condition of (ii) in loc.cit follows from 
 Proposition 1.9 in Chap. II of \cite{dS}.

\subhead{Acknowledgement}\endsubhead The possibility 
to  construct systems like this was suggested by Kazuya Kato.

\no Kenichiro Kimura

\no Institute of Mathematics

\no  University of Tsukuba

\no Tsukuba, Ibaraki, 305 

\no Japan 

\no email: kimurak\@math.tsukuba.ac.jp

\enddocument